\documentclass[12pt,a4paper]{article} 

\usepackage{amsmath,amssymb,graphicx}
\usepackage[french,english]{babel}

\newtheorem{theoreme}{Theorem}[section]

\newtheorem{define}[theoreme]{Definition}
\newtheorem{propal}[theoreme]{Proposition}
\newtheorem{corolle}[theoreme]{Corollary}

\newenvironment{proof2}
    {\begin{list}{{\bf Proof}}
                 {\setlength{\labelwidth}{0cm}       
                  \setlength{\leftmargin}{0em}}    
                 \item}
    {$\blacksquare$\end{list}}
\newenvironment{remarque}{\par \addvspace{\bigskipamount} 
                               \refstepcounter{theoreme}
                               {\bf Remark \thetheoreme ~:~}}
                         {\par \addvspace{\bigskipamount}}

\def \disp {\displaystyle}

\def \N {{\mathbb N}} 
\def \R {{\mathbb R}}
\def \C {{\cal C}}

\def \eps {{\varepsilon}}

\def \A {{\cal A}}
\def \D {{\cal D}}
\def \Nor {{\cal N}}
\def \E {{{\cal E}}}
\def \S {{\cal S}}

\def \un {{\ensuremath{\textrm{\textnormal{1}}\hspace{-0.7ex}
                       \textrm{\textnormal{I}}}}}

\def \bx {{\bf x}} \def \nx {{\breve{x}}}
\def \by {{\bf y}} \def \ny {{\breve{y}}}
\def \bz {{\bf z}} \def \nz {{\breve{z}}} 
\def \bv {{\bf v}} \def \nv {{\breve{v}}}
\def \bX {{\bf X}} \def \nX {{\breve{X}}}
\def \bW {{\bf W}} \def \nW {{\breve{W}}}
\def \bn {{\bf n}} \def \nn {{\breve{n}}} 
\def \bl {{\bf l}}
\def \bL {{\bf L}} 
\def \bsigma {{{\boldsymbol\sigma}}} \def \nsigma {{\breve{\sigma}}}
\def \bb {{\bf b}} \def \nb {{\breve{b}}} 

\def \rm {r_-}
\def \rp {r_+}

\def \phi {{\varphi}}

\author{
Myriam \textsc{Fradon} \\
Universit\'e de Lille1
}

\date{28th October 2009}

\title{Brownian Dynamics of Globules}

\begin{document}

\maketitle

\centerline{\textbf{Abstract}}
We prove the existence and uniqueness of a strong solution of a stochastic differential equation with normal reflection representing the random motion of finitely many globules. Each globule is a sphere with time-dependent random radius and a center moving according to a diffusion process. The spheres are hard, hence non-intersecting, which induces in the equation a reflection term with a local (collision-)time. 
A smooth interaction is considered too and, in the particular case of a gradient system, the reversible measure of the dynamics is given. 
In the proofs, we analyze geometrical properties of the boundary of the set in which the process takes its values, in particular the so-called  Uniform Exterior Sphere and Uniform Normal Cone properties.
These techniques extend to other hard core models of objects with a time-dependent random characteristic: we present here an application to the random motion of a chain-like molecule.\\


\section{Introduction}

Since the pioneering work of Skorokhod \cite{Skorokhod1et2}, many authors have investigated the question of the existence and uniqueness of a solution for reflected stochastic differential equations in a domain. It has first been solved for half-spaces, then for convex domains (see \cite{Tanaka}). Lions and Sznitman \cite{LionsSznitman} proved the existence of a solution in so-called \emph{admissible sets} and Saisho \cite{SaishoSolEDS} extended these results to domains satisfying only the \emph{Uniform Exterior Sphere} and the \emph{Uniform Normal Cone} conditions (see definitions in Section \ref{resultSaisho}).
These results were applied to prove the existence and uniqueness of Brownian dynamics for hard spheres in \cite{SaishoTanakaBrownianBalls}, \cite{FR2}, \cite{FR3}, or for systems of mutually reflecting molecules (see \cite{SaishoMolecules2}).

We are interested here in dynamics of finitely many objects having not only a random position but also another random time-dependent geometrical characteristic like the radius for spheres or the length of the bonds in a molecule. We will prove the existence and uniqueness of random dynamics for two elaborated models, using methods which are refinements of Saisho's techniques, analyzing fine geometrical properties of the boundary of the domain on which the motion is reflected.

More precisely~:

We first introduce a globules model representing a finite system of non-intersecting spheres which centers undergo diffusions and which radii vary according to other diffusions constrained to stay between a maximum and a minimum value. 
The spheres might be cells, or particles, or soap bubbles floating on the surface of a liquid (2-dimensional motion) or in the air (3-dimensional motion). 
The behavior of the globules is quite intuitive~:
two globules collide when the distance between their random centers is equal to the sum of their random radii, and the collision has as effect that their centers  move away from one another and their sizes decrease. The associated stochastic differential equation $(\E_g)$ (see Section \ref{Sect_globules_model}) includes several reflection terms, each of them corresponding to a constraint in the definition of the set of allowed globules configurations~: constraints of non-intersection between each couple of spheres, constraints on the radii to stay between fixed minimum and maximum values. We proves that this equation has a unique strong solution and give in some special case a time-reversible initial distribution (theorems \ref{thexistglobules} and \ref{threversglobules}).

We also consider a model for linear molecules, such as alkanes (carbon chains) or polymers~:
each atom moves like a diffusion, the lengths of the bonds between neighbour atoms vary between a minimum and maximum value which evolve according to a reflected diffusion. 
This corresponds to a SDE $(\E_c)$ reflected on the boundary of the set of all allowed chains.
Here also, we prove the existence and uniqueness of the solution of $(\E_c)$ (see theorem \ref{thexistchenilles}) with similar methods as in the globule case.

The rest of the paper is organized as follows~: In Section \ref{SectGeomCriterium}, we present a new general criterion for a multiple constraint domain to satisfy the Uniform Exterior Sphere and Uniform Normal Cone conditions (see proposition \ref{propcritere}). These geometrical assumptions on the boundary induce existence and uniqueness of the reflected diffusion on this domain.
We also obtain a disintegration of the local time along the different constraint directions .
Section \ref{SectPreuves} is devoted to the proofs of the theorems announced in Section \ref{Sect_two_models}.

For the sake of shortness, we restricted ourselves to these two examples, though dynamics for Brownian systems evolving under multiple constraints may be found in other situations (the results in section \ref{SectGeomCriterium} are given in a general frame for easier adaptation to other examples).

\section{Two hard core models}
\label{Sect_two_models}
\subsection{Globules model}
\label{Sect_globules_model}

We want to construct a model for interacting globules. Each globule is spherical with random radius oscillating between a minimum and a maximum value. Its center is a point in $\R^d$, $d \ge 2$. The number $n$ of globules is fixed. Globules configurations will be denoted by
$$
\bx=(x_1,\nx_1,\ldots,x_n,\nx_n) 
\quad\text{ with }\quad x_1,\ldots,x_n\in\R^d \quad\text{ and }\quad \nx_1,\ldots,\nx_n\in \R
$$
where $x_i$ is the center of the $i^\text{th}$ globule and $\nx_i$ is its radius.
An allowed globules configuration is a configuration $\bx$ satisfying
$$
\forall i \quad \rm\le \nx_i \le\rp \quad\text{ and }\quad
\forall i\neq j \quad |x_i-x_j|\ge \nx_i+\nx_j
$$ 
So, in an allowed configuration, spheres do not intersect and their radii are bounded from below by the minimum value $\rm>0$ and bounded from above by the maximum value $\rp>\rm$. 
In this paper, the symbol $|\cdot|$ denotes the Euclidean norm on $\R^d$ or $\R^{(d+1)n}$ (or some other Euclidean space, depending on the context).

Let $\A_g$ be the set of allowed globules configurations~:
$$
\A_g=\left\{ \bx \in \R^{(d+1)n},~\forall i \quad \rm\le \nx_i \le\rp \text{ and } \forall i\neq j \quad |x_i-x_j|\ge \nx_i+\nx_j 
     \right\}
$$
The random motion of reflecting spheres with fluctuating radii is represented by the following stochastic differential equation~:
$$
(\E_g)
\left\{
\begin{array}{l}
 \disp
  X_i(t) = X_i(0) + \int_0^t \sigma_i(\bX(s)) dW_i(s) + \int_0^t b_i(\bX(s)) ds      
           + \sum_{j=1}^n \int_0^t \frac{X_i(s)-X_j(s)}{\nX_i(s)+\nX_j(s)} dL_{ij}(s) \\
 \disp
 \nX_i(t) = \nX_i(0) + \int_0^t \nsigma_i(\bX(s)) d\nW_i(s) 
            + \int_0^t \nb_i(\bX(s)) ds - \sum_{j=1}^n L_{ij}(t) - L^+_i(t) + L^-_i(t)
\end{array}
\right.
$$
In this equation, $\bX(s)$ is the vector $(X_i(s),\nX_i(s))_{1 \le i\le n}$.
The $W_i$'s are independents $\R^d$-valued Brownian motions and the $\nW_i$'s are independents one-dimensional Brownian motions, also independent from the $W_i$'s.
The diffusion coefficients $\sigma_i$ and $\nsigma_i$, and the drift coefficients $b_i$ and $\nb_i$ are functions defined on $\A_g$, with values in the $d \times d$ matrices for $\sigma_i$, values in $\R^d$ for $b_i$, and values in $\R$ for $\nsigma_i$ and $\nb_i$.
To make things simpler with the summation indices, we let $ L_{ii} \equiv 0 $.

A solution of equation $(\E_g)$ is a continuous $\A_g$-valued process $\{ \bX(t) , t \ge 0 \}$ satisfying equation $(\E_g)$ for some family  of local times $L_{ij}$, $L^+_i$, $L^-_i$ such that for each $i,j$~:
$$
(\E'_g)
\left\{
\begin{array}{l}
\disp
L_{ij} \equiv L_{ji}, \quad  
L_{ij}(t) = \int_0^t \un_{|X_i(s)-X_j(s)|=\nX_i(s)+\nX_j(s)} ~dL_{ij}(s)~, \\
\disp
L^+_i(t) = \int_0^t \un_{\nX_i(s)=\rp} ~dL^+_i(s) \quad\textrm{ and }\quad L^-_i(t) = \int_0^t \un_{\nX_i(s)=\rm} ~dL^-_i(s) 
\end{array}
\right.
$$ 
\begin{remarque} \label{remtempslocal}
Here, and in the sequel, the expression {\bf local time} stands for~: non-decreasing adapted continuous process, which starts from $0$ and has bounded variations on each finite interval.
\end{remarque}
Equation $(\E_g)$ has an intuitive meaning~:
\begin{itemize}
\item
the positions and radii of the spheres are Brownian;
\item
when a globule becomes too big, it's deflated~: $\nX_i$ decreases by $-dL^+_i$ when $\nX_i=\rp$;
\item
when a globule becomes too small, it's inflated~: $\nX_i$ increases by $+dL^-_i$ when $\nX_i=\rm$;
\item
when two globules bump into each other (i.e. $|X_i-X_j|=\nX_i+\nX_j$), they are deflated and move away from each other~: 
$\nX_i$ decreases by $-dL_{ij}$ and $X_i$ is given an impulsion in the direction $\frac{X_i-X_j}{\nX_i+\nX_j}$ with an amplitude $dL_{ij}$. 
\end{itemize}

In the case of an hard core interaction between spheres with a \emph{fixed} radius, the existence of solutions for the corresponding SDE has been proved in \cite{SaishoTanakaBrownianBalls}. However, the condition $|x_i-x_j|\ge \nx_i+\nx_j$ is not equivalent to $|(x_i,\nx_i)-(x_j,\nx_j)|\ge c$ for some real number $c$. Hence the above model is not a classical hard sphere model in $\R^d$ or $\R^{d+1}$. 

\begin{theoreme}
\label{thexistglobules}
Assume that the diffusion coefficients $\sigma_i$ and $\nsigma_i$ and the drift coefficients $b_i$ and $\nb_i$ are bounded and Lipschitz continuous on $\A_g$ (for $1 \le i \le n$).
Then equation $(\E_g)$ has a unique strong solution.
\end{theoreme}

\begin{remarque}
"Strong uniqueness of the solution" here stands for strong uniqueness (in the sense of \cite{IkedaWatanabe} chap.IV def.1.6) of the process $\bX$, and, as a consequence, of the reflection term.
This does not imply strong uniqueness for the local times $L_{ij}$, $L^+_i$, $L^-_i$ unless several collisions at the same time with linearly dependant collision directions does not occurs a.s. (see the proof of corollary \ref{corolSaisho} for details).
\end{remarque}

The first part of the next theorem is a corollary of the previous one. The second part describes the equilibrium states of systems of interacting globules. Here, and through this paper, $d\bx$ denotes the Lebesgue measure.

\begin{theoreme}
\label{threversglobules}
The diffusion representing the motion of $n$ globules submitted to a smooth interaction $\phi$ exists as soon as the interaction potential $\phi$ is an even $\C^2$ function on $\R^d$ with bounded derivatives. It is the unique strong solution of the equation~:
$$
(\E^{\phi}_g)
\left\{
\begin{array}{l}
 \disp
 X_i(t) = X_i(0) + W_i(t) -\frac{1}{2} \int_0^t \sum_{j=1}^n \nabla\phi(X_i(s)-X_j(s)) ds      
         + \sum_{j=1}^n \int_0^t \frac{X_i(s)-X_j(s)}{\nX_i(s)+\nX_j(s)} dL_{ij}(s) \\
 \disp
 \nX_i(t) = \nX_i(0) + \nW_i(t) - \sum_{j=1}^n L_{ij}(t) - L^+_i(t) + L^-_i(t) \\
 \text{with } \bX~~ \A_g\text{-valued continuous process and } (L_{ij},L^+_i,L^-_i)_{0 \le i,j \le n} \text{ satisfying conditions } (\E'_g)
\end{array}
\right.
$$
Moreover, if $Z=\int_{\A_g} e^{-\sum_{1 \le i<j \le n} \phi(x_i-x_j)} d\bx <+\infty$ then the Probability measure $\mu$ defined by $d\mu(\bx)=\frac{1}{Z} \un_{\A_g}(\bx) e^{-\sum_{1 \le i<j \le n} \phi(\bx_i-\bx_j)} d\bx$ is time-reversible for this diffusion, i.e. $(\bX(T-t))_{t\in[0,T]}$ has the same distribution as $(\bX(t))_{t\in[0,T]}$ for each positive $T$ when $\bX(0) \sim \mu$.
\end{theoreme}

These theorems are proved in section \ref{PreuvesGlobules}. The proofs rely on previous results from Y.Saisho and H. Tanaka (see section \ref{resultSaisho}) and on an inheritance criterion for geometrical properties which is given in section \ref{criteredheredite}.

\subsection{Linear molecule model}

Another example of hard core interaction between particles with another spatial characteristic beside position is the following simple model for a linear molecule. In this model, we study chains of particles having a fixed number of links with variable length. More precisely, a configuration is a vector
\[
\bx=(x_1,\ldots,x_n,\nx_-,\nx_+)
\quad\text{ with }\quad x_1,\ldots,x_n\in\R^d \quad\text{ and }\quad \nx_-,\nx_+\in \R
\]
$x_i$ and $x_{i+1}$ are the ends of the $i^\text{th}$ link and $\nx_->0$ (resp. $\nx_+ > \nx_-$) is the minimum (resp. maximum) allowed length of the links for the chain. The number $n$ of particles in the chain is at least equal to $2$, they are moving in $\R^d$, $d \ge 2$.
So the set of allowed configurations is~:
$$
\A_c=\left\{ \bx \in \R^{dn+2},~ \rm \le \nx_- \le \nx_+ \le \rp \text{ and } 
             \forall i\in\{1,\ldots,n-1\} \quad \nx_- \le |x_i-x_{i+1}| \le \nx_+ \right\}
$$
We want to construct a model for the random motion of such chains, as the $\A_c$-valued solution of the following stochastic differential equation~:
$$
(\E_c)
\left\{
\begin{array}{l}
 \disp
  X_i(t) = X_i(0) + \int_0^t \sigma_i(\bX(s)) dW_i(s) + \int_0^t b_i(\bX(s)) ds 
                  + \int_0^t \frac{X_i-X_{i+1}}{\nX_-}(s) dL^-_i(s)                                     \\
\phantom{X_i(t) =}\disp
           + \int_0^t \frac{X_i-X_{i-1}}{\nX_-}(s) dL^-_{i-1}(s)
           - \int_0^t \frac{X_i-X_{i+1}}{\nX_+}(s) dL^+_i(s) - \int_0^t \frac{X_i-X_{i-1}}{\nX_+}(s) dL^+_{i-1}(s) \\
 \disp
 \nX_-(t) = \nX_-(0) + \int_0^t \sigma_-(\bX(s)) d\nW_-(s) + \int_0^t b_-(\bX(s)) ds
                     - \sum_{i=1}^{n-1} L^-_i(t) + L_-(t) - L_=(t)  \\
 \disp
 \nX_+(t) = \nX_-(0) + \int_0^t \sigma_+(\bX(s)) d\nW_+(s) + \int_0^t b_+(\bX(s)) ds
                     + \sum_{i=1}^{n-1} L^+_{i}(t) - L_+(t) + L_=(t)
\end{array}
\right.
$$
As in the previous model, the $W_i$'s are independents $\R^d$-valued Brownian motions and the $\nW_-$ and $\nW_+$ are independents one-dimensional Brownian motions, also independent from the $W_i$'s.
A solution of equation $(\E_c)$ is an $\A_c$-valued continuous process $\{ \bX(t), t \ge 0 \}$ satisfying the equation for some family of 
local times $L^-_i$, $L^+_i$, $L_-$, $L_+$, $L_=$ such that for each $t \in \R^+ $ and each $i$~:
$$
\begin{array}{l}
\disp
L^-_i(t) = \int_0^t \un_{|X_i(s)-X_{i+1}(s)|=\nX_-(s)} dL^-_i(s)~,\quad 
L^+_i(t) = \int_0^t \un_{|X_i(s)-X_{i+1}(s)|=\nX_+(s)} dL^+_i(s)     \\
\disp
L_-(t) = \int_0^t \un_{\nX_-(s)=\rm} dL_-(s) ,\quad L_+(t) = \int_0^t \un_{\nX_+(s)=\rp} dL_+(s) \quad\text{and}\quad
L_=(t) = \int_0^t \un_{\nX_-(s)=\nX_+(s)} dL_=(s)
\end{array}
$$ 
Equation $(\E_c)$ looks more complicated than $(\E_g)$ because it contains a larger number of local times, but it is very simple on an intuitive level~: both ends of each link are Brownian, links that are too short ($|X_i-X_{i+1}|=\nX_-$) tend to become longer (reflection term in direction $X_i-X_{i+1}$ for $X_i$ and in the opposite direction for $X_{i+1}$) and also tend to diminish the lower bound $\nX_-$ (negative reflection term in the equation of $\nX_-$). 
Symmetrically, links that are too large ($|X_i-X_{i+1}|=\nX_+$) are both becoming shorter (reflection term in direction $X_{i+1}-X_i$ for $X_i$) and enlarging the lower bound $\nX_+$ (positive reflection term in its equation).
Moreover, $\nX_-$ increases (by $L_-$) if it reaches its lower limit $\rm$ and $\nX_+$ decreases (by $L_+$) if it reaches its lower limit $\rp$. And the lower bound decreases and the upper bound increases (by $L_=$) when they are equal, so as to fulfill the condition $\nx_- \le \nx_+$.

\begin{theoreme}
\label{thexistchenilles}
If the $\sigma_i$'s, $\sigma_-$ and $\sigma_+$ and the $b_i$'s, $b_-$ and $b_+$ are bounded and Lipschitz continuous on $\A_c$, then equation $(\E_c)$ has a unique strong solution.

Moreover, assume that the $\sigma_i$'s, $\sigma_-$ and $\sigma_+$ are equal to the identity matrix, $b_-$ and $b_+$ vanish, and $b_i(\bx)=-\frac{1}{2}\sum_{j=1}^n \nabla\phi(x_i-x_j)$ for some even $\C^2$ function  $\phi$ on $\R^d$ with bounded derivatives, satisfying $Z=\int_{\A_c} e^{-\sum_{1 \le i<j \le n} \phi(x_i-x_j)} d\bx <+\infty$. Then the solution with initial distribution 
$d\mu(\bx)=\frac{1}{Z} \un_{\A_c}(\bx) e^{-\sum_{1 \le i<j \le n} \phi(\bx_i-\bx_j)} d\bx$ is time-reversible.
\end{theoreme}
See section \ref{PreuvesChenilles} for the proof of this theorem.

\section{Geometrical criteria for the existence of reflected dynamics} \label{SectGeomCriterium}

\subsection{Uniform Exterior Sphere and Uniform Normal Cone properties} \label{resultSaisho}

In order to solve the previous stochastic differential equations, we will use theorems of Saisho and Tanaka extending some previous results of  Lions and Sznitman \cite{LionsSznitman}. 

To begin with, we need geometrical conditions on subset boundaries. The subsets we are interested in are sets of allowed configurations. But for the time being, we just consider any subset $\D$ in $\R^m$($m\ge 2$) which is the closure of an open connected set with non-zero (possibly infinite) volume. $d\bx$ is the Lebesgue measure on $\R^m$, $|\cdot|$ denotes the Euclidean norm as before, and we set
$$
\S^m=\{ \bx\in\R^m,~~|\bx|=1 \}
$$
We define the set of all (inward) normal vectors at point $\bx$ on the boundary $\partial\D$ as~:
\[
\Nor^{\D}_\bx=\bigcup_{\alpha>0} \Nor^{\D}_{\bx,\alpha} \quad\text{ where }\quad 
\Nor^{\D}_{\bx,\alpha}=\{\bn\in\S^m,~~ \mathring{B}(\bx-\alpha\bn,\alpha) \cap \D = \emptyset \}
\]
Here $\mathring{B}(\bx,r)$ is the open ball with radius $r$ and center $\bx$.
Note that
\[
\mathring{B}(\bx-\alpha\bn,\alpha) \cap \D = \emptyset
\quad\quad\Longleftrightarrow\quad\quad
\forall \by\in\D~~~ (\by-\bx).\bn+\frac{1}{2\alpha}|\by-\bx|^2 \ge 0
\]
\begin{define}
If there exists a constant $\alpha>0$ such that $\Nor^{\D}_\bx=\Nor^{\D}_{\bx,\alpha} \neq \emptyset$ for each $\bx\in\partial\D$, we say that $\D$ has the {\bf Uniform Exterior Sphere} property (with constant $\alpha$) and we write~: $\D\in UES(\alpha)$.
\end{define}
$UES(\alpha)$ means that a sphere of radius $\alpha$ rolling on the outside of $\D$ can touch each point of $\partial\D$.
This property is weaker than the convexity property (which corresponds to $UES(\infty)$) but still ensures the existence of a local  projection function similar to the projection on convex sets. 

\begin{define}
We say that $\D$ has the {\bf Uniform Normal Cone} property with constants $\beta$, $\delta$, and we write $\D\in UNC(\beta,\delta)$,
if for some $\beta\in [0,1[$ and $\delta>0$, for each $\bx\in\partial\D$, there exists $\bl_\bx \in \S^m$ such that for every $\by\in\partial\D$
\[
|\by-\bx| \le \delta \quad\Longrightarrow\quad \forall \bn\in\Nor^{\D}_\by \quad \bn.\bl_\bx \ge \sqrt{1-\beta^2}
\]
\end{define}

Let us consider the reflected stochastic differential equation~:
\begin{equation}
\bX(t)=\bX(0) + \int_0^t \bsigma(\bX(s)) d\bW(s) + \int_0^t \bb(\bX(s)) ds + \int_0^t \bn(s)d\bL(s) \label{eqSaisho}
\end{equation}
where $(\bW(t))_t$ is a $m$-dimensional Brownian motion. 
A solution of (\ref{eqSaisho}) is a $(\bX,\bn,\bL)$ with $\bX$ a $\D$-valued adapted continuous process, $\bn(s)\in\Nor^{\D}_{\bX(s)}$ when $\bX(s)\in\partial\D$ and $\bL$ a local time such that
\[
\bL(\cdot)=\int_0^{\cdot} \un_{\bX(s)\in\partial\D} d\bL(s)
\]
The following result is a consequence of \cite{SaishoSolEDS} and \cite{SaishoTanakaSymmetry}~:
\begin{theoreme}\label{thSaisho}
Assume $\D\in UES(\alpha)$ and $\D\in UNC(\beta,\delta)$ for some positive $\alpha$, $\delta$ and some $\beta\in[0,1[$. 
If $\sigma:\D\longrightarrow\R^{m^2}$ and $b:\D\longrightarrow\R^m$ are bounded Lipschitz continuous functions, then  (\ref{eqSaisho}) has a unique strong solution. 
Moreover, if $\bsigma$ is the identity matrix and $\disp \bb=-\frac{1}{2} \nabla\Phi$ with $\Phi$ a $\C^2$ function on $\R^m$ with bounded derivatives satisfying $Z=\int_{\D} e^{-\Phi(\bx)} d\bx <+\infty$, then the solution with initial distribution $d\mu(\bx)=\frac{1}{Z} \un_{\D}(\bx) e^{-\Phi(\bx)} d\bx$ is time-reversible.
\end{theoreme}

\begin{proof2} {\bf of theorem \ref{thSaisho}} \\
The existence of a unique strong solution of (\ref{eqSaisho}) is proved in \cite{SaishoSolEDS} theorem 5.1.

Let us consider the special case of the gradient system~: $\bb=-\frac{1}{2} \nabla\Phi$ and $\bsigma$ identity matrix.
For fixed $T>0$ and $\bx\in\D$, let $P_{\bx}$ be the distribution of the solution $(\bX,\bn,\bL)$ of (\ref{eqSaisho}) starting from $\bx$. The process $(W(t))_{t\in[0,T]}$ defined as
$$
\bW(t)=\bX(t)-\bx +\frac{1}{2} \int_0^t \nabla\Phi(\bX(s)) ds - \int_0^t \bn(s)d\bL(s)
$$
is a Brownian motion with respect to the probability measure $P_{\bx}$.
Since $\nabla\Phi$ is bounded, thanks to Girsanov theorem, the process $\disp \tilde{\bW}(t)=\bW(t)-\frac{1}{2} \int_0^t \nabla\Phi(\bX(s)) ds$ is a Brownian motion with respect to the probability measure $\tilde{P}_{\bx}$ defined by $\disp \frac{d\tilde{P}_{\bx}}{dP_{\bx}}=M_T$ where
$$
M_t=\exp\left( \frac{1}{2}\int_0^t \nabla\Phi(\bX(s)).d\bW(s) - \frac{1}{8} \int_0^t |\nabla\Phi(\bX(s))|^2 ds \right)
$$
From theorem 1 of \cite{SaishoTanakaSymmetry}, it is known that Lebesgue measure on $\D$ is time-reversible for the solution 
of $\disp \bX(t)=\bX(0) + \tilde{\bW}(t) + \int_0^t \bn(s)d\bL(s)$, that is, the measure $\disp \tilde{P}_{d\bx}=\int_{\D} \tilde{P}_{\bx} d\bx$ is invariant under time reversal on $[0;T]$ for any positive $T$. Using It\^o's formula to compute~:
{\small$$
\int_0^t \nabla\Phi(\bX(s)).d\bW(s)
=\Phi(\bX(T))-\Phi(\bX(0)) +\frac{1}{2}\int_0^T |\nabla\Phi(\bX(s))|^2-\Delta\Phi(\bX(s))~ds -\int_0^T \nabla\Phi(\bX(s)).\bn(s) d\bL(s) 
$$}
we notice that the density of the probability measure $\disp P_{\mu}=\int_{\D} P_{\bx} \mu(d\bx)$ with respect to the measure $\tilde{P}_{d\bx}$ is equal to~:
{\small
$$
\frac{e^{-\Phi(\bX(0))}}{Z} \exp\left( \frac{\Phi(\bX(0))-\Phi(\bX(T))}{2}  
                                       +\int_0^T \frac{1}{4} \Delta\Phi(\bX(s)) -\frac{1}{8}|\nabla\Phi(\bX(s))|^2 ds 
                                       +\frac{1}{2} \int_0^T \nabla\Phi(\bX(s)).\bn(s) d\bL(s) \right)
$$}
This expression does not change under time reversal, i.e. if $(\bX,\bn,\bL)$ is replaced by $(\bX(T-\cdot),\bn(T-\cdot),\bL(T)-\bL(T-\cdot)$. Thus the time reversibility of  $\tilde{P}_{d\bx}$ implies the time reversibility of $P_{\mu}$.
\end{proof2}

\subsection{The special case of multiple constraints}\label{criteredheredite}

The sets we are interested in are sets of configurations satisfying multiple constraints. They are intersections of several sets, each of them defined by a single constraint. So we need a sufficient condition on the $\D_i$'s for Uniform Exterior Sphere and Uniform Normal Cone properties to hold on $\D=\cap_{i=1}^p \D_i$.

In \cite{SaishoMolecules2}, Saisho gives a sufficient condition for the intersection $\D=\cap_{i=1}^p \D_i$ to inherit both properties when each set $\D_i$ has these properties. However, in order to prove the existence of a solution in the case of our globules model, we have to consider the intersection of the sets~:
\[
\D_{ij}=\left\{ \bx \in \R^{dn+n},~ |x_i-x_j|\ge \nx_i+\nx_j \right\}
\]
Uniform Exterior Sphere property does not hold for $\D_{ij}$'s, so Saisho's inheritance criterion does not work here. However, $\D_{ij}$ satisfies an Exterior Sphere condition restricted to $\A_g$ in some sense, and we shall check that this is enough. Similarly, the Uniform Normal Cone property does not hold for $\D_{ij}$, but a restricted version holds, and proves sufficient for our needs. 

So we present here a UES and UNC criterion with weaker assumptions. We keep Saisho's smoothness assumption, which holds for many interesting models, and restrict the UES and UNC assumptions on the $\D_i$'s to the set $\D$. Instead of the conditions $(B_0)$ and $(C)$ in \cite{SaishoMolecules2}, we introduce the compatibility  assumption (iv) which is more convenient because it is easier to check a property at each point $\bx\in\partial\D$ for finitely many vectors normal to the $\partial\D_i$'s, than on a neighborhood of each $\bx$ for the whole (infinite) set of vectors normal to $\partial\D$. 

The proof of this criterion is postponed to the end of this section.

\begin{propal} {\bf (Inheritance criterion for UES and UNC conditions)} \\
\label{propcritere}
For $\D\subset\R^m$ equal to the intersection $\disp \D=\bigcap_{i=1}^p \D_i$, assume that~:
\begin{itemize}
\item[(i)]
The sets $\D_i$ are closures of domains with non-zero volumes and boundaries at least $\C^2$ in $\D$~: this implies the existence of a unique unit normal vector $\bn_i(\bx)$ at each point $\bx\in \D\cap\partial\D_i$.
\item[(ii)]
Each set $\D_i$ has the Uniform Exterior Sphere property restricted to $\D$, i.e. \\
$ \disp
\quad \exists \alpha_i>0 , \quad \forall \bx\in \D\cap\partial\D_i \quad \mathring{B}(\bx-\alpha_i\bn_i(\bx),\alpha_i) \cap \D_i = \emptyset
$
\item[(iii)]
Each set $\D_i$ has the Uniform Normal Cone property restricted to $\D$, i.e. \\
for some $\beta_i\in [0,1[$ and $\delta_i>0$ and for each $\bx\in\D\cap\partial\D_i$ there is a unit vector $\bl^i_\bx$ s.t.
$\quad \forall \by\in\D\cap\partial\D_i \cap B(\bx,\delta_i) \quad \bn_i(\by).\bl^i_\bx \ge \sqrt{1-\beta_i^2}$
\item[(iv)] (compatibility assumption) There exists $\beta_0>\sqrt{2\max_{1 \le i \le p} \beta_i}$ satisfying\\
$\disp \forall\bx\in\partial\D,~~~ \exists \bl_\bx^0 \in \S^m,
\quad \forall i \text{ s.t. } \bx\in\partial\D_i \quad \bl_\bx^0.\bn_i(\bx) \ge \beta_0$
\end{itemize}
Under the above assumptions, 
$\D\in UES(\alpha)$ and $\D\in UNC(\beta,\delta)$ hold with $\alpha=\beta_0 \min_{1 \le i \le p} \alpha_i$,
$\delta=\min_{1 \le i \le p} \delta_i/2$ and $\beta=\sqrt{1-(\beta_0-2 \max_{1 \le i \le p} \beta_i / \beta_0)^2}$.
Moreover, the vectors normal to the boundary $\partial\D$ are convex combinations
of the vectors normal to the boundaries $\partial\D_i$~:
\[
\forall \bx\in\partial\D \quad 
\Nor^\D_\bx=\left\{ \bn\in\S^m,~~~ \bn=\sum_{\partial\D_i\ni\bx} c_i \bn_i(\bx) \text{ with each } c_i \ge 0 \right\}
\]
\end{propal}

\begin{remarque}\label{remsommeci}
Thanks to the compatibility assumption $(iv)$, $\bn=\sum_{\partial\D_i\ni\bx} c_i \bn_i(\bx)$ with non-negatives $c_i$'s and $|\bn|=1$ implies that
$$
\sum_{\partial\D_i\ni\bx} c_i \le \sum_{\partial\D_i\ni\bx} c_i \frac{\bn_i(\bx).\bl_\bx^0}{\beta_0} = \frac{\bn.\bl_\bx^0}{\beta_0} \le \frac{1}{\beta_0}
$$
so that the last equality in proposition \ref{propcritere} can be rewritten as
$$
\Nor^\D_\bx=\left\{ \bn\in\S^m,~~~ \bn=\sum_{\partial\D_i\ni\bx} c_i \bn_i(\bx)
                    \text{ with } \forall i ~~ c_i \ge 0 \text{ and } \sum_{\partial\D_i\ni\bx} c_i \le \frac{1}{\beta_0}\right\}
$$
\end{remarque}

\begin{corolle} {\bf of theorem \ref{thSaisho}} \\
\label{corolSaisho}
If $\disp \D=\bigcap_{i=1}^p \D_i$ satisfies assumptions $(i)\cdots(iv)$ and if $\bsigma$ and $\bb$ are bounded Lipschitz continuous functions, then
\begin{equation}
\label{eqSaishoLi}
\bX(t)=\bX(0) +\int_0^t \bsigma(\bX(s)) d\bW(s) +\int_0^t \bb(\bX(s)) ds + \sum_{i=1}^p \int_0^t \bn_i(\bX(s)) dL_i(s)
\end{equation}
has a unique strong solution with local times $L_i$ satisfying
$\disp
L_i(\cdot) = \int_0^\cdot \un_{\partial\D_i}(\bX(s)) ~dL_i(s) $
\end{corolle}

\begin{proof2} {\bf of corollary \ref{corolSaisho}} \\
Thanks to proposition \ref{propcritere}, $\D$ satisfies the assumptions of theorem \ref{thSaisho}.
Thus equation (\ref{eqSaishoLi}) has a unique strong solution $\bX$ with local time $\bL$ and reflection direction $\bn$. Using proposition \ref{propcritere} again, there are (non-unique) coefficients $c_i(\omega,s)\in [0,\frac{1}{\beta_0}]$ for each normal vector in the reflection term to be written as a convex combination~:
\begin{equation}\label{egaliteci}
\bn(\omega,s)=\sum_{\partial\D_i\ni\bX(\omega,s)} c_i(\omega,s) \bn_i(\bX(\omega,s)) 
\end{equation}
Let us prove that there exists a measurable choice of the $c_i$'s~:

The map $\bn$ (resp. $\bn_i(\bX)$) is only defined for $(\omega,s)$ such that $\bX(\omega,s)\in\partial\D$ (resp. $\bX(\omega,s)\in\partial\D_i$). 
We extend these maps by zero to obtain measurable maps on $\Omega\times[0,T]$ (for an arbitrary positive $T$).
Note that equality (\ref{egaliteci}) holds for the extended maps too.

For each $(\omega,s)$, we define the map $f_{\omega,s}$ on $\R^p$ by $f_{\omega,s}(c)=|\sum_{i=1}^p c_i \bn_i(\bX(\omega,s)) - \bn(\omega,s)|$.
For a positive integer $k$, let $R_k=\{0,\frac{1}{k},\frac{2}{k},\ldots,\lfloor \frac{k}{\beta_0} \rfloor \}^p$ denote the $\frac{1}{k}$-lattice on $[0,\frac{1}{\beta_0}]^p$ endowed with lexicographic order. The smaller point in $R_k$ for which $f_{\omega,s}$ reaches its minimum value is
$$
c^{(k)}(\omega,s)=\sum_{c\in R_k} c \prod_{c'\in R_k, c'\neq c} \left( \un_{f_{\omega,s}(c')>f_{\omega,s}(c)} 
                                                             + \un_{f_{\omega,s}(c')=f_{\omega,s}}(c) \un_{c'>c} \right)
$$
$c^{(k)}$ is a measurable map on $\Omega\times[0,T]$ and (\ref{egaliteci}) implies that $|f_{\omega,s}(c^{(k)}(\omega,s))| \le \frac{p}{k}$.
Taking (coordinate after coordinate) the limsup of the sequence of $(c^{(k)})_k$, we obtain a measurable process $c^{(\infty)}$ satisfying (\ref{egaliteci}).

Finally, let $\disp L_i = \int_0^{\cdot} \un_{\partial\D_i}(\bX(s)) c_i(s)~d\bL(s)$. 
$L_i$ has bounded variations on each $[0,T]$ since the $c_i$'s are bounded, and it is a local time in the sense of remark \ref{remtempslocal}.

Here, strong uniqueness holds for process $\bX$ and for the reflection term $\disp \sum_{i=1}^p \int_0^t \bn_i(\bX(s)) dL_i(s)$.
The uniqueness of this term does not imply uniqueness of the $L_i$'s, because the $c_i$'s are not unique. 
\end{proof2}

\begin{proof2} {\bf of proposition \ref{propcritere}} \\
Let $\bx\in\partial\D$. Since $\D=\bigcap_{i=1}^p \D_i$, the set $\{i \text{ s.t. } \bx\in\partial\D_i \}$ is not empty.
Since $\D_i$ satisfies $UES(\alpha_i)$ restricted to $\D$, for each $\by\in\D_i$~~~$(\by-\bx).\bn_i(\bx)+\frac{1}{2\alpha_i}|\by-\bx|^2 \ge 0$, and consequently~:
\[
\forall i \text{ s.t. } \partial\D_i\ni\bx \quad \forall \by\in\D \quad (\by-\bx).\bn_i(\bx) + \frac{1}{2\min_{\partial\D_j\ni\bx}\alpha_j} |\by-\bx|^2 \ge 0
\]
Summing over $i$ for non-negative $c_i$'s such that $\sum_{\partial\D_i\ni\bx} c_i~\le\frac{1}{\beta_0}$ we obtain~:
\[
\forall \by\in\D \quad 
(\sum_{\partial\D_j\ni\bx} c_i \bn_i(\bx)).(\by-\bx)+\frac{1}{2 \beta_0 \min_{\partial\D_i\ni\bx}\alpha_j} |\by-\bx|^2 \ge 0
\]
which implies that $\sum_{\partial\D_j\ni\bx} c_i \bn_i(\bx)$ belongs to the set $\Nor_{\bx,\alpha}$ of normal vectors on the boundary of $\D$, for $\alpha=\beta_0 \min_{\partial\D_j\ni\bx}\alpha_j$. Thanks to remark \ref{remsommeci}, this proves the inclusion
\[
\Nor'_\bx:=\left\{ \bn\in\S^m,~~~ \bn=\sum_{\partial\D_i\ni\bx} c_i \bn_i(\bx) \text{ with } c_i \ge 0 \right\}
~~~\subset~~~ \Nor_{\bx,\alpha}
\]
Let us prove the converse inclusion.
We first remark that the smoothness of the boundary $\partial\D_i$ at point $\bx$ implies that for any $\eps>0$~:
\[
\exists \delta^i_\eps(\bx)>0 \quad\text{ s.t. }\quad 
\{ \bx+\bz,~~ |\bz|\le\delta^i_\eps(\bx) \text{ and } \bz.\bn_i(\bx) \ge\eps|\bz| \} \subset \D_i
\]
Consequently, for $\disp N_\eps=\bigcap_{\partial\D_i\ni\bx} \{\bz,~\bn_i(\bx).\bz \ge \eps|\bz| \}$ we have~:
$
\{ \bx+\bz,~~ \bz\in N_\eps \text{ and } |\bz|\le\min_{i}\delta^i_\eps \} \subset \D $

By definition, for each $\bn\in\Nor^\D_\bx$, there exist $\alpha_\bn>0$ such that $\forall \by\in\D~~ (\by-\bx).\bn+\frac{1}{2\alpha_\bn}|\by-\bx|^2 \ge 0$, hence for $\bz\in N_\eps$ and $\lambda>0$ small enough~:
\[
\lambda \bz.\bn+\frac{\lambda^2}{2\alpha_\bn}|\bz|^2 \ge 0
\]
For this to hold even with $\lambda$ going to zero, $\bz.\bn$ has to be non-negative. So we obtain~:
\[
\forall \bn\in \Nor^\D_\bx \quad \forall \eps>0 \quad -\bn \in N_\eps^*
\]
where $N_\eps^*=\{\bv,~~\forall \bz\in N_\eps ~\bv.\bz \le 0 \}$ is the dual cone of the convex cone $N_\eps$.
As proved in Fenchel \cite{Fenchel} (see also \cite{SaishoTanakaBrownianBalls}), the dual of a finite intersection of convex cones is the set of all limits of linear combinations of their dual cones, in particular~:
\[
N_\eps^* =\overline{ \sum_{\partial\D_i\ni\bx} \{\bz,~\bn_i(\bx).\bz \ge \eps|\bz| \}^* }
         =\sum_{\partial\D_i\ni\bx} \{ \bv,~ -\bn_i(\bx).\bv \ge \sqrt{1-\eps^2} |\bv| \}
\]
For $k\in \N$ large enough, since $-\bn \in N_{\frac{1}{k}}^*$, there exist unit vectors $\bn_{i,k}$ and non-negative numbers $c_{i,k}$ such that~:
\[
\bn=\sum_{\partial\D_i\ni\bx} c_{i,k} \bn_{i,k} \quad\quad \text{ and } \quad\quad
\text{ for } \partial\D_i\ni\bx \quad \bn_i(\bx).\bn_{i,k} \ge \sqrt{1-\frac{1}{k^2}} 
\] 
When $k$ tends to infinity, $\bn_{i,k}$ tends to $\bn_i(\bx)$, and for $k$ large enough $\bn_{i,k}.\bl^0_\bx \ge \frac{\beta_0}{2}$ thus~:
\[
1 \ge \bn.\bl^0_\bx \ge \sum_{\partial\D_i\ni\bx} c_{i,k} \bn_{i,k}.\bl^0_\bx \ge \frac{\beta_0}{2} \sum_{\partial\D_i\ni\bx} c_{i,k}
\]
Thus the sequences $(c_{i,k})_k$ are bounded, which implies the existence of convergent subsequences. Their limits $c_{i,\infty} \ge 0$ satisfy~:
\[
\bn=\sum_{\partial\D_i\ni\bx} c_{i,\infty} \bn_i(\bx)
\]
This completes the proof of $\Nor^\D_\bx \subset \Nor'_\bx$.
We already proved that $\Nor'_\bx \subset \Nor^\D_{\bx,\alpha}$ with $\alpha=\beta_0 \min_{\partial\D_j\ni\bx}\alpha_j$, so we obtain 
$\Nor^\D_\bx = \Nor'_\bx = \Nor^\D_{\bx,\alpha}$ for each $\bx\in\partial\D$.
As a consequence, $\D\in UES(\beta_0 \min_{1 \le j \le p} \alpha_j)$.

Let us now prove that $\D\in UNC(\beta,\delta)$.
The Uniform Normal Cone property restricted to $\D$ holds for $\D_i$, with constant $\beta_i<\frac{\beta_0^2}{2} \le \frac{1}{2}$.
That is, for $\bx\in\D\cap\partial\D_i$ there exist a unit vector $\bl^i_\bx$ which satisfies $\bl^i_\bx.\bn_i(\by) \ge \sqrt{1-\beta_i^2}$ for each $\by\in\D\cap\partial\D_i$ such that $|\bx-\by|\le\delta_i$. \\
If $\bl^i_\bx=\bn_i(\bx)$, this implies that $|\bn_i(\bx)-\bn_i(\by)|^2 \le 2\beta_i^2$. \\
If $\bl^i_\bx \neq \bn_i(\bx)$, we use the Gram-Schmidt orthogonalization process for a sequence of vectors with $\bl^i_\bx$ and $\bn_i(\bx)$ as first vectors, then compute $\bn_i(\bx).\bn_i(\by)$ in the resulting orthonormal basis $(\bl^i_\bx,{\bf e}_2,{\bf e}_3,\ldots,{\bf e}_m)$~:
\[
\bn_i(\bx).\bn_i(\by)=(\bl^i_\bx.\bn_i(\by))(\bl^i_\bx.\bn_i(\bx))+(\bn_i(\by).{\bf e}_2)\sqrt{1-(\bl^i_\bx.\bn_i(\bx))^2}
\]
Note that $|\bn_i(\by).{\bf e}_2|\le\beta_i$ because $\bl^i_\bx.\bn_i(\by) \ge \sqrt{1-\beta_i^2}$ and $|\bn_i(\by)|=1$, thus~:
\[
\bn_i(\bx).\bn_i(\by) \ge \Big(\sqrt{1-\beta_i^2}\Big)^2-\beta_i^2=1-2\beta_i^2 
\]
This implies that $|\bn_i(\bx)-\bn_i(\by)|^2 \le 4\beta_i^2$.\\
So in both cases~: $\disp |\bn_i(\bx)-\bn_i(\by)| \le 2\beta_i$ as soon as $|\bx-\by|\le\delta_i$ for $\bx,\by\in\D\cap\partial\D_i$.

Let us now fix $\bx\in\partial\D$ and $\delta=\min_{1 \le i \le p} \delta_i/2$.
We then choose $\bx'\in \partial\D \cap B(\bx,\delta)$ such that $\{i \text{ s.t. } \bx'\in\partial\D_i \} \supset \{i \text{ s.t. } \by\in\partial\D_i \}$ for each $\by\in \partial\D \cap B(\bx,\delta)$ and we let $\bl=\bl^0_{\bx'}$. 
To complete the proof of $\D\in UNC(\beta,\delta)$, we only have to prove that $\bn.\bl$ is uniformly bounded from below for $\bn\in\Nor^\D_\by$ with $\by\in \partial\D \cap B(\bx,\delta)$.

We already know that each $\bn\in\Nor^\D_\by$ is a convex sum of elements of the $\Nor^{\D_i}_\by$~:
$\bn=\sum_{\partial\D_i\ni\by} c_i \bn_i(\by) $.
The coefficients $c_i$ are non-negative, their sum is not smaller than $1$ because $\bn$ is a unit vector, and is not larger than $\frac{1}{\beta_0}$ thanks to $(iv)$.
So the vector $\bn'=\sum_{\partial\D_i\ni\by} c_i \bn_i(\bx')$ satisfies~:
\[
\bn'.\bl=\sum_{\partial\D_i\ni\bx'} c_i \bn_i(\bx').\bl \ge \beta_0
\quad\text{ and }\quad
|\bn'-\bn| \le \sum_{\partial\D_i\ni\by} c_i |\bn_i(\bx')-\bn_i(\by)| 
\le 2 \sum_{\partial\D_i\ni\by} c_i \beta_i \le 2 \max_{1 \le i \le p} \frac{\beta_i}{\beta_0}
\]
Consequently~:
$\bn.\bl \ge \bn'.\bl-|\bn'-\bn| \ge \beta_0-2 \max_{1 \le i \le p} \frac{\beta_i}{\beta_0} >0$.
\end{proof2}

\section{Existence of dynamics for globules and linear molecule models } \label{SectPreuves}

\subsection{Globules model} \label{PreuvesGlobules}

Let us prove that the globules model satisfies the assumptions in proposition \ref{propcritere}.
The set of allowed configurations is~:
\[
\begin{array}{c} \disp
\disp \A_g=(\bigcap_{1 \le i<j \le n} \D_{ij}) \cap (\bigcap_{1 \le i \le n} \D_{i+}) \cap (\bigcap_{1 \le i  \le n} \D_{i-}) \\~\\
\text{ where } \quad \D_{ij}=\left\{ \bx \in \R^{dn+n},~ |x_i-x_j|\ge \nx_i+\nx_j \right\}                                    \\~\\
                     \D_{i+}=\left\{ \bx \in \R^{dn+n},~ \nx_i \le \rp \right\} \quad 
                     \D_{i-}=\left\{ \bx \in \R^{dn+n},~ \nx_i \ge \rm \right\}
\end{array}
\]
The $\D_{ij}$ have smooth boundaries on $\A_g$ and the characterization of normal vectors is easy~: at point $\bx$ satisfying $|x_i-x_j|=\nx_i+\nx_j>0$, the unique unit inward normal vector $\bn_{ij}(\bx)=\bn$ is given by~:   
\[
n_i=\frac{x_i-x_j}{2(\nx_i+\nx_j)} \quad \nn_i=-\frac{1}{2} \quad n_j=\frac{x_j-x_i}{2(\nx_i+\nx_j)} \quad \nn_j=-\frac{1}{2}
\quad \text{(every other component vanishes)}
\]
On the boundary of the half-space $\D_{i+}$, at point $\bx$ such that $\nx_i=\rp$, the unique unit normal vector $\bn$ has only one non-zero component~: $\nn_i=-1$. Similarly, $\D_{i-}$ is a half-space, $\bx$ belongs to its boundary if $\nx_i=\rm$, and the unique unit normal vector $\bn$ at this point has $\nn_i=1$ as its only non-zero component.
These vectors do not depend on $\bx$ and will be denoted by $\bn_{i+}$, $\bn_{i-}$ instead of $\bn_{i+}(\bx)$, $\bn_{i-}(\bx)$.
\begin{propal}\label{PropGeomBulles}
$\A_g$ satisfies properties $\disp UES\left( \frac{\rm^2}{2\rp n\sqrt{n}} \right)$ and $\disp UNC\left( \sqrt{1-\frac{\rm^2}{2^6\rp^2 n^3}}~,~\frac{\rm^5}{2^{14}\rp^4 n^6} \right)$.\\
Moreover, the vectors normal to the boundary $\partial\A_g$ are convex combinations of the vectors normal to the boundaries $\partial\D_{ij}$, $\partial\D_{i+}$, $\partial\D_{i-}$, that is, for every $\bx$ in $\partial\A_g$~:
\[
\Nor^{\A_g}_\bx=\left\{ \bn\in\S^{dn+n},~~
                    \bn=\sum_{ \partial\D_{ij}\ni\bx} c_{ij} \bn_{ij}(\bx) 
                    +\sum_{ \partial\D_{i+}\ni\bx } c_{i+} \bn_{i+}
                    +\sum_{ \partial\D_{i-}\ni\bx } c_{i-} \bn_{i-}
                    \text{ with } c_{ij},c_{i+},c_{i-} \ge 0 \right\}
\]
\end{propal}

\begin{proof2} {\bf of proposition \ref{PropGeomBulles}}

We have to check that the assumption of proposition \ref{propcritere} are satisfied for the set $\A_g$.

Since the $\D_{i+}$ and $\D_{i-}$ are half-spaces, the Uniform Exterior Sphere property holds for them with any positive constant (formally $\alpha_{i+}=\alpha_{i-}=+\infty$). For the same reason, the Uniform Normal Cone property holds for $\D_{i+}$ with any constants $\beta_{i+}$ and $\delta_{i+}$, and with $\bl^{i+}_\bx$ equal to the normal vector $\bn_{i+}$. This also holds for the sets $\D_{i-}$ with any $\beta_{i-}$ and $\delta_{i-}$, and with $\bl^{i-}_\bx=\bn_{i+}$.

Let us consider $\bx\in\A_g$ such that $\bx\in\partial\D_{ij}$, i.e. $|x_i-x_j|=\nx_i+\nx_j$. 
By definition of $\bn_{ij}(\bx)$, for $\by=\bx-(\nx_i+\nx_j)\bn_{ij}(\bx)$, one has~:
\[
\begin{array}{c} \disp
y_i=x_i-(\nx_i+\nx_j)\frac{x_i-x_j}{2(\nx_i+\nx_j)}=\frac{x_i+x_j}{2}=x_j-(\nx_i+\nx_j)\frac{x_j-x_i}{2(\nx_i+\nx_j)}=y_j \\
\ny_i+\ny_j=\nx_i-(\nx_i+\nx_j)(-\frac{1}{2}) +\nx_j-(\nx_i+\nx_j)(-\frac{1}{2})=2(\nx_i+\nx_j)
\end{array}
\]
For $\bz\in\mathring{B}(0,\nx_i+\nx_j)$~:
\[
|(y_i+z_i)-(y_j+z_j)|-(\ny_i+\nz_i+\ny_j+\nz_j) \le |z_i|+|z_j|-2(\nx_i+\nx_j)+|\nz_i|+|\nz_j| \le 2|\bz|-2(\nx_i+\nx_j) <0
\]
thus $\by+\bz \in\D_{ij}^c$. This proves that $\mathring{B}(\by,\nx_i+\nx_j) \subset \D_{ij}^c$, hence $\Nor^{\D_{ij}}_\bx= \Nor^{\D_{ij}}_{\bx,\nx_i+\nx_j}$. The general Uniform Exterior Sphere property does not hold, but the property restricted to $\A_g$ holds for $\D_{ij}$ with constant $\alpha_{ij}=2\rm$ because $\nx_i+\nx_j \ge 2\rm$ for $\bx\in\A_g$.

For $\bx\in\A_g$ such that $\bx\in\partial\D_{ij}$, let us define $\bl^{ij}_\bx=\bn_{ij}(\bx)$.
For another configuration $\by\in\partial\D_{ij}$~:
\[
\bl^{ij}_\bx.\bn_{ij}(\by)
=\frac{x_i-x_j}{2(\nx_i+\nx_j)}.\frac{y_i-y_j}{2(\ny_i+\ny_j)}+\frac{x_j-x_i}{2(\nx_i+\nx_j)}.\frac{y_j-y_i}{2(\ny_i+\ny_j)}
+(-\frac{1}{2})^2+(-\frac{1}{2})^2
=\frac{(x_i-x_j).(y_i-y_j)}{2(\nx_i+\nx_j)(\ny_i+\ny_j)}+\frac{1}{2}
\]
Since $(x_i-x_j).(y_i-y_j) \ge |x_i-x_j|^2-|x_i-x_j||y_i-x_i-y_j+x_j| \ge |x_i-x_j|^2-\sqrt{2}|x_i-x_j||\bx-\by|$
and $\ny_i+\ny_j \le (\nx_i+\nx_j)+\sqrt{2}|\bx-\by|$, this leads to~:
\[
\bl^{ij}_\bx.\bn_{ij}(\by) \ge \frac{|x_i-x_j|-\sqrt{2}|\bx-\by|}{2(\nx_i+\nx_j+\sqrt{2}|\bx-\by|)}+\frac{1}{2}
=\frac{\nx_i+\nx_j}{\nx_i+\nx_j+\sqrt{2}|\bx-\by|}
\]
Consequently $\bl^{ij}_\bx.\bn_{ij}(\by) \ge \sqrt{1-\beta_{ij}^2}$ as soon as 
$|\bx-\by| \le \frac{\nx_i+\nx_j}{\sqrt{2}} \left( \frac{1}{\sqrt{1-\beta_{ij}^2}}-1 \right)$.
This proves that $\D_{ij}$ has the Uniform Normal Cone property restricted to $\A_g$ with any constant $\beta_{ij}\in ]0,1[$ and with $\delta_{ij}=\sqrt{2} \rm \left( \frac{1}{\sqrt{1-\beta_{ij}^2}}-1 \right)>0$. In particular, the property holds for any $\beta_{ij}$ with $\delta_{ij}=\rm \beta_{ij}^2/\sqrt{2}$. 

For $\bx\in\partial\A_g$, let us construct a unit vector $\bl_\bx^0$ satisfying assumption \emph{(iv)} in proposition \ref{propcritere}.
We first have to define clusters of colliding globules~:
\[
C_\bx(i)=\left\{ j \text{ s.t. } 
                 |x_i-x_{j_1}|=\nx_i+\nx_{j_1},~ |x_{j_1}-x_{j_2}|=\nx_{j_1}+\nx_{j_2},\ldots,
                 |x_{j_k}-x_j|=\nx_{j_k}+\nx_j \text{ for some } j_1,\ldots,j_k 
         \right\}
\] 
and their centers of gravity $\disp x'_i=\frac{1}{\sharp C_\bx(i)} \sum_{j\in C_\bx(i)} x_j$. 
Let $\bl_\bx^0=\frac{\bv}{|\bv|}$ with for each $i$~:
\[
v_i= x_i-x'_i \quad\text{ and }\quad \nv_i=\frac{\rm}{\rp-\rm} (\frac{\rp+\rm}{2}-\nx_i)
\] 
We need an upper bound on the norm of $\bv$. Since $\bx\in\A_g$, each $\nx_i$ is larger than $\rm$ and smaller than $\rp$, thus $|\frac{\rp+\rm}{2}-\nx_i| \le \frac{\rp-\rm}{2}$. Moreover, if $j$ belongs to the cluster $C_\bx(i)$ around the $i^\text{th}$ globule, the distance between $x_i$ and $x_j$ is at most $2(n-1)\rp$, thus $|v_i|\le 2(n-1)\rp$~:
\[
|\bv|^2=\sum_{i=1}^n |v_i|^2 + (\frac{\rm}{\rp-\rm})^2 \sum_{i=1}^n (\frac{\rp+\rm}{2}-\nx_i)^2
\le 4n(n-1)^2 \rp^2 + n\frac{\rm^2}{4}  < 4 n^3 \rp^2
\]
We also need lower bounds on the scalar products of $\bv$ with normal vectors on the boundaries of the $\D_{i+}$, $\D_{i-}$, $\D_{ij}$. If $\nx_i=\rp$, then $\nv_i=\frac{\rm}{\rp-\rm} \frac{\rm-\rp}{2}$ thus $\bv.\bn_{i+}=\frac{\rm}{2}$. Similarly, if $\nx_i=\rm$, then $\bv.\bn_{i-}=\frac{\rm}{2}$. If $|x_i-x_j|=\nx_i+\nx_j$ then $x'_i=x'_j$ thus~:
\[
\begin{array}{l}
\bv.\bn_{ij}(\bx)=(x_i-x'_i).\frac{x_i-x_j}{2(\nx_i+\nx_j)}+(x_j-x'_j).\frac{x_j-x_i}{2(\nx_i+\nx_j)}
                    -\frac{1}{2}\frac{\rm}{\rp-\rm} (\frac{\rp+\rm}{2}-\nx_i)
                    -\frac{1}{2}\frac{\rm}{\rp-\rm} (\frac{\rp+\rm}{2}-\nx_j)                             \\
\phantom{\bv.\bn_{ij}(\bx)}
=(x_i-x_j).\frac{x_i-x_j}{2(\nx_i+\nx_j)}-\frac{1}{2}\frac{\rm}{\rp-\rm}(\rp+\rm-\nx_i-\nx_j)             \\
\phantom{\bv.\bn_{ij}(\bx)}
=\frac{\rp}{2(\rp-\rm)}(\nx_i+\nx_j)-\frac{\rm(\rp+\rm)}{2(\rp-\rm)} \ge \frac{\rm}{2}
\end{array}
\]
because $\nx_i+\nx_j \ge 2\rm$. This proves that assumption \emph{(iv)} in proposition \ref{propcritere} is satisfied with 
$\disp \beta_0 = \frac{\rm}{2\sqrt{4 n^3 \rp^2}}=\frac{\rm}{4\rp n\sqrt{n}}$.
The inequality $\beta_0>\sqrt{2\max(\max_i \beta_{i+},\max_i \beta_{i-},\max_{i,j} \beta_{ij})}$ holds as soon as
$\beta_{ij}<\beta_0^2/2$. As seen above, the Uniform Normal Cone property restricted to $\A_g$ holds for $\D_{ij}$ with constants $\beta_{ij}=\beta_0^2/4$ and $\delta_{ij} =\frac{\rm\beta_0^4}{16\sqrt{2}}$.
So, thanks to proposition \ref{propcritere}, we obtain that $\A_g$ has the Uniform Exterior Sphere property with constant 
$\alpha_{\A_g}=2\rm\beta_0=\frac{\rm^2}{2\rp n\sqrt{n}}$ and the Uniform Normal Cone property with constants 
$\delta_{\A_g}=\frac{\rm\beta_0^4}{32\sqrt{2}}=\frac{\rm^5}{2^{13}\sqrt{2}\rp^4 n^6}$ and
$\beta_{\A_g}=\sqrt{1-(\beta_0/2)^2}=\sqrt{1-\frac{\rm^2}{64\rp^2 n^3}}$
\end{proof2}

\begin{proof2} {\bf of theorems \ref{thexistglobules} and \ref{threversglobules}}

As seen in the proof of proposition \ref{PropGeomBulles}, the set $\A_g$ of allowed globules configurations satisfies the assumptions of corollary \ref{corolSaisho}.
If the diffusion coefficients $\sigma_i$ and $\nsigma_i$ and the drift coefficients $b_i$ and $\nb_i$ are bounded and Lipschitz continuous on $\A_g$ (for $1 \le i \le n$), then the functions
$$
\bsigma(\bx)=\left[\begin{array}{ccccc} \sigma_1(\bx)&    0         & \cdots & \cdots      &   0          \\
                                            0        &\nsigma_1(\bx)& \ddots &             & \vdots       \\
                                          \vdots     & \ddots       & \ddots & \ddots      & \vdots       \\
                                          \vdots     &              & \ddots &\sigma_n(\bx)&  0           \\
                                            0        & \cdots       & \cdots &   0         &\nsigma_n(\bx) 
                    \end{array}\right]
\text{ and }
\bb(\bx)=\left[\begin{array}{c} b_1(\bx)   \\
                                \nb_1(\bx) \\
                                \vdots     \\
                                b_n(\bx)   \\
                                \nb_n(\bx) 
                                \end{array}\right]
$$
are bounded and Lipschitz continuous as well. Thus, for $m=n(d+1)$ and $\D=\A_g$, equation (\ref{eqSaisho}) has a unique strong solution $\bX$, and there exists a decomposition of the reflection term $\int_0^{\cdot} \bn(s)d\bL(s)$ so that it is equal to 
$$
\left( 
\sum_{j=1}^n \int_0^{\cdot} \frac{X_i-X_j}{2(\nX_i+\nX_j)}(s) dL_{ij}(s)~,~ 
-\sum_{j=1}^n \int_0^{\cdot} \frac{1}{2} dL_{ij}(s) -\int_0^{\cdot} dL_{i+}(s) +\int_0^{\cdot} dL_{i-}(s) 
\right)_{1 \le i \le n}
$$
with $\disp L_{ij} = \int_0^{\cdot} \un_{|X_i(s)-X_j(s)|=\nX_i(s)+\nX_j(s)} c_{ij}(s)~d\bL(s)$,\quad
$\disp L^+_i = \int_0^{\cdot} \un_{\nX_i(s)=\rm} c_{i+}(s)~d\bL(s)$ and 
$\disp L^-_i = \int_0^{\cdot} \un_{\nX_i(s)=\rp} c_{i-}(s)~d\bL(s)$ 
for some measurable choice of the $c_{ij},c_{i+},c_{i-}$.
For convenience, the local time $\frac{1}{2}L_{ij}$ has been used in equation $(\E_g)$, and denoted by $L_{ij}$ again.

For an even $\C^2$ function $\phi$ on $\R^{d+1}$ with bounded derivatives, 
the function $\Phi(\bx)=\sum_{1 \le i<j \le n} \phi(\bx_i-\bx_j)$ is $\C^2$ with bounded derivatives
and the drift function in equation $(\E^{\phi}_g)$ is equal to $-\frac{1}{2} \nabla\Phi$. 
Thus $(\E^{\phi}_g)$ has a unique strong solution.
Moreover, if $Z=\int_{\A_g} e^{-\Phi(\bx)} d\bx <+\infty$, theorem \ref{thSaisho} implies the time-reversibility of the solution with initial distribution $d\mu(\bx)=\frac{1}{Z} \un_{\A_g}(\bx) e^{-\Phi(\bx)} d\bx$.
\end{proof2}

\subsection{Linear molecule model} \label{PreuvesChenilles}

The set of allowed chains configurations is the intersection of $2n+1$ sets~:
\[
\begin{array}{c} 
\disp \A_c= (\bigcap_{i=1}^{n-1} \D_{i-}) \bigcap~ (\bigcap_{i=1}^{n-1} \D_{i+}) \bigcap~ \D_- \bigcap~ \D_+ \bigcap~ \D_= \\~\\
\disp\text{ where } \quad \D_{i-} =\left\{ \bx \in \R^{dn+2},~ |x_i-x_{i+1}|\ge \nx_- \right\} \quad
                          \D_{i+} =\left\{ \bx \in \R^{dn+2},~ |x_i-x_{i+1}|\le \nx_+ \right\}                             \\~\\
\disp                     \D_- =\left\{ \bx \in \R^{dn+2},~ \nx_- \ge \rm \right\} \quad
                          \D_+ =\left\{ \bx \in \R^{dn+2},~ \nx_+ \le \rp \right\} \quad
                          \D_= =\left\{ \bx \in \R^{dn+2},~ \nx_- \le \nx_+ \right\}
\end{array}
\]
The boundaries of these sets are smooth, and simple derivation computations give the unique unit vector normal to each boundary at each point of $\partial\A_c$ (we are not interested in other points). Actually, $\D_-$, $\D_+$ and $\D_=$ are half-spaces, which makes the computations and checking of UES and UNC properties very simple.
\begin{itemize}
\item
The vector $\bn_-=(0,\ldots,0,1,0)$ is normal to the boundary $\partial\D_-$ at each point $\bx\in\partial\D_-$
\item
Similarly, $\bn_+=(0,\ldots,0,0,-1)$ is normal to the boundary $\partial\D_+$ at each point $\bx\in\partial\D_+$
\item
The vector normal to the boundary $\partial\D_=$ at each point $\bx\in\partial\D_=$ is 
$\bn_= =(0,\ldots,0,\frac{-1}{\sqrt{2}},\frac{1}{\sqrt{2}})$
\item
At every $\bx\in\A_c\cap\partial\D_{i-}$ (for $1 \le i \le n-1$), the  vector normal to $\partial\D_{i-}$ is $\bn=\bn_{i-}(\bx)$ defined by~:
\[
n_i=\frac{x_i-x_{i+1}}{\nx_- \sqrt{3}} \quad n_{i+1}=\frac{x_{i+1}-x_i}{\nx_- \sqrt{3}} \quad \nn_-=-\frac{1}{\sqrt{3}}
\quad \text{(other components equal zero)}
\]
\item
At every $\bx\in\A_c\cap\partial\D_{i+}$, the  vector normal to $\partial\D_{i+}$ is $\bn=\bn_{i+}(\bx)$ defined by~:
\[
n_i=\frac{x_{i+1}-x_i}{\nx_+ \sqrt{3}} \quad n_{i+1}=\frac{x_i-x_{i+1}}{\nx_+ \sqrt{3}} \quad \nn_+=\frac{1}{\sqrt{3}}
\quad \text{(other components equal zero)}
\]
\end{itemize}
The proof of theorem \ref{thexistchenilles} is similar to the proof of theorems \ref{thexistglobules} and \ref{threversglobules}, and will be omitted. It relies on the following proposition and uses theorem \ref{thSaisho} and corollary \ref{corolSaisho} to obtain the existence, uniqueness, and reversibility of the solution of $(\E_c)$. As in the proof of theorem \ref{thexistglobules}, the local times are multiplied by a suitable constant to provide a more convenient expression.

So, in order to prove theorem \ref{thexistchenilles}, we only have to check that $\A_c$ satisfies the assumptions of proposition \ref{propcritere}, i.e. that the following proposition holds~:
\begin{propal}\label{PropGeomCaterpillar}
$\A_c$ satisfies properties $\disp UES\left(\frac{\rm^2}{2\rp n \sqrt{2n}}\right)$ and 
$\disp UNC\left(\sqrt{1-\frac{\rm^2}{24\rp^2 n^3}},\frac{\rm^5 \sqrt{3}}{2^{12} \rp^4 n^6}\right)$. \\
Moreover, for each $\bx$ in $\partial\A_c$, the set $\Nor^{\A_c}_\bx$ of normal vectors is equal to~:
{\small\[
\left\{ 
\begin{array}{r} \disp
\bn\in\S^{dn+2},
\bn=\sum_{\partial\D_{i-}\ni\bx} c_{i-} \bn_{i-}(\bx) +\sum_{\partial\D_{i+}\ni\bx} c_{i+} \bn_{i+}(\bx) 
    + \un_{\partial\D_-}(\bx) c_- \bn_- + \un_{\partial\D_+}(\bx)c_+ \bn_+ + \un_{\partial\D_=}(\bx) c_= \bn_= \\
\text{ with } c_{i-},c_{i+},c_-,c_+,c_= \ge 0 
\end{array}
\right\}
\]}
\end{propal}

\begin{proof2} {\bf of proposition \ref{PropGeomCaterpillar}}

Let us check that the set $\A_c$ satisfies the assumptions of proposition \ref{propcritere}.

The Uniform Exterior Sphere property and the Uniform Normal Cone property hold for $\D_-$, $\D_+$ and $\D_=$, with any positive $\alpha$ and $\delta$ and any $\beta$ in $[0,1[$, because these sets are half-spaces.

For $\bx\in\A_c\cap\partial\D_{i-}$, for $\alpha=\frac{\rm\sqrt{3}}{2}$, and for $\by=\bx-\alpha\bn_{i-}(\bx)$~:
\[
y_i=x_i-\alpha\frac{x_i-x_{i+1}}{\nx_- \sqrt{3}} \quad y_{i+1}=x_{i+1}-\alpha\frac{x_{i+1}-x_i}{\nx_- \sqrt{3}} \quad 
\ny_-=\nx_- +\frac{\alpha}{\sqrt{3}} 
\]
so that each $\bz\in\mathring{B}(0,\alpha)$ satisfies
$|(y_i+z_i)-(y_{i+1}+z_{i+1})|=\left| x_i-x_{i+1}- \rm\frac{x_i-x_{i+1}}{\nx_-} + z_i-z_{i+1} \right|$ 
and that $|x_i-x_{i+1}|=\nx_-$ implies 
\[
|(y_i+z_i)-(y_{i+1}+z_{i+1})|-(\ny_- + \nz_-) 
\le \nx_-(1-\frac{\rm}{\nx_-}) +|z_i-z_{i+1}|-\nx_- -\frac{\rm}{2}+ \nz_-
=-\alpha\sqrt{3}+|z_i-z_{i+1}|+ \nz_-
\]
We know that $|z_i-z_{i+1}|+ \nz_- \le \sqrt{3}|\bz| < \sqrt{3}\alpha$, thus the above quantity is negative and $\by+\bz \not\in \D_{i-}$.
Consequently, the Uniform Exterior Sphere property holds for $\D_{i-}$ with constant $\alpha_{i-}=\frac{\rm\sqrt{3}}{2}$.
It also holds for $\D_{i+}$, with any positive constant, because this set is  convex (the midpoint of two points in $\D_{i+}$ obviously belongs to $\D_{i+}$).
In order to prove that $\D_{i-}$ has the Uniform Normal Cone property restricted to $\A_c$, we fix two chains $\bx,\by\in\A_c\cap\partial\D_{i-}$ and compute~:
\[
\bn_{i-}(\bx).\bn_{i-}(\by)=\frac{2}{3}\frac{x_i-x_{i+1}}{\nx_-}.\frac{y_i-y_{i+1}}{\ny_-}+\frac{1}{3}
=\frac{2|x_i-x_{i+1}|^2}{3 \nx_-^2}+\frac{1}{3}
  +\frac{2}{3}\frac{x_i-x_{i+1}}{\nx_-}.\left( \frac{y_i-y_{i+1}}{\ny_-}-\frac{x_i-x_{i+1}}{\nx_-} \right)
\]
Since $|x_i-x_{i+1}|=\nx_-$, we obtain~:
$\disp \bn_{i-}(\bx).\bn_{i-}(\by) \ge 1-\frac{2}{3}\left| \frac{y_i-y_{i+1}}{\ny_-}-\frac{x_i-x_{i+1}}{\nx_-} \right|$.\\
Then $|\nx_- (y_i-y_{i+1})-\ny_- (x_i-x_{i+1})| \le \nx_-|y_i-y_{i+1}-x_i+x_{i+1}|+|\nx_- - \ny_-||x_i-x_{i+1}|$ leads to~:
\[
\bn_{i-}(\bx).\bn_{i-}(\by) \ge 1-\frac{2}{3 \ny_-}\left( |y_i-x_i|+|y_{i+1}-x_{i+1}|+|\nx_- - \ny_-|\right)
\ge 1-\frac{2|\by-\bx|}{\sqrt{3} \ny_-}
\ge 1-\frac{2|\by-\bx|}{\sqrt{3} \rm}
\]
Thus $\bn_{i-}(\bx).\bn_{i-}(\by) \ge \sqrt{1-\beta_{i-}^2}$ as soon as $|\by-\bx| \le \frac{\rm\sqrt{3}}{2}(1-\sqrt{1-\beta_{i-}^2})$, which is implied by $|\by-\bx| \le \frac{\rm\sqrt{3}}{4}\beta_{i-}^2$.
As a consequence, the set $\D_{i-}$ has the Uniform Normal Cone property restricted to $\A_c$ with any constant $\beta_{i-}\in[0,1[$ and the corresponding constant $\delta_{i-}=\frac{\rm\sqrt{3}}{4}\beta_{i-}^2$. 

Note that for $\bx,\by\in\A_c\cap\partial\D_{i+}$~:
$ \disp\quad \bn_{i+}(\bx).\bn_{i+}(\by)=\frac{2}{3}\frac{x_i-x_{i+1}}{\nx_+}.\frac{y_i-y_{i+1}}{\ny_+}+\frac{1}{3}$ \\
so that the same computation gives that $\D_{i+}$ has the Uniform Normal Cone property restricted to $\A_c$ with any constants $\beta_{i+}\in[0,1[$ and $\delta_{i+}=\frac{\rm\sqrt{3}}{4}\beta_{i+}^2$. 

In order to check the compatibility assumption, let us fix a point $\bx\in\partial\A_c$ and construct a suitable vector $\bv$ such that $\bl_\bx^0=\frac{\bv}{|\bv|}$.
We first define the "middle point" $x'$ of the chain~: $x'=x_{(n+1)/2}$ if the chain contains an odd number of particles, and
$x'=\frac{x_{n/2}+x_{n/2+1}}{2}$ if $n$ is an even number. We also define $\nx'=\frac{\nx_+ + \nx_-}{2}$.
We then construct the $v_i$'s in an incremental way, starting from the middle and going up and down to both ends~:
\begin{itemize}
\item
If $n$ is odd, we choose $v_{(n+1)/2}=0$.
\item
In the case of an even $n$, we choose $(v_{n/2},v_{n/2+1})=(x_{n/2}-x',x_{n/2+1}-x')$ if $|x_{n/2}-x_{n/2+1}|<\nx'$ and 
$(v_{n/2},v_{n/2+1})=(x_{n/2+1}-x',x_{n/2}-x')$ if $|x_{n/2}-x_{n/2+1}|>\nx'$. \\
In the critical case $|x_{n/2}-x_{n/2+1}|=\nx'$, our choice depends on the value of $\nx'$~: we let $(v_{n/2},v_{n/2+1})=(x_{n/2}-x',x_{n/2+1}-x')$ if $\nx'<\rp$ and $(v_{n/2},v_{n/2+1})=(x_{n/2+1}-x',x_{n/2}-x')$ if $\nx_+=\nx_-=\rp$.
\item
The other $v_i$'s are chosen incrementally so as to fulfill the same condition~:
\[
v_i-v_{i+1}=\left\{ \begin{array}{l} 
                    x_i-x_{i+1} ~~~\text{ if } |x_i-x_{i+1}|<\nx' \text{,~~ or } |x_i-x_{i+1}|=\nx' \text{ and } \nx'<\rp \\
                    x_{i+1}-x_i ~~~\text{ if } |x_i-x_{i+1}|>\nx' \text{,~~ or } |x_i-x_{i+1}|=\nx' \text{ and } \nx'=\rp
                    \end{array}    \right.
\]
\end{itemize}
Note that $|v_i|\le (n-1)\rp$ for each $i$. For our choice of $\bv$ to be complete, we also define~:
\begin{itemize}
\item
$\disp \nv_- =\frac{\rm}{2},~~ \nv_+ =\frac{-\rm}{2} $ \quad if $\nx_- < \nx_+$
\item
$\disp \nv_- =\frac{\rm}{2},~~ \nv_+ =\frac{3}{2}\rp $ \quad if $\nx_- = \nx_+ < \rp$
\item
$\disp \nv_- =\frac{-3}{2}\rp,~~ \nv_+ =\frac{-\rp}{2} $ \quad if $\nx_- = \nx_+ = \rp$
\end{itemize}
As in the proof of proposition \ref{PropGeomBulles}, we need an upper bound on the norm of $\bv$~:
\[
|\bv|^2=\sum_{i=1}^n |v_i|^2 + \nv_-^2 + \nv_+^2 \le n(\frac{n-1}{2}\rp)^2 + \frac{\rp^2}{4} + \frac{9\rp^2}{4} \le \frac{n^3 \rp^2}{2}
\]
We also have to prove that the scalar products of $\bv$ with normal vectors on the boundaries are uniformly bounded from below. 
Note that $\bv.\bn_- =\nv_-$,\quad $\bv.\bn_+ =-\nv_+$,\quad $\bv.\bn_= = \frac{\nv_+ - \nv_-}{\sqrt{2}}$,\quad 
$\disp \bv.\bn_{i-}(\bx)=\frac{x_i-x_{i+1}}{\nx_- \sqrt{3}}.(v_i-v_{i+1}) - \frac{\nv_-}{\sqrt{3}}$\quad  and \quad 
$\disp \bv.\bn_{i+}(\bx)=\frac{x_i-x_{i+1}}{\nx_+ \sqrt{3}}.(v_{i+1}-v_i) + \frac{\nv_+}{\sqrt{3}}$.

In the case where $\nx_- < \nx_+$, the choices made on $\bv$ lead to~:
\begin{itemize}
\item
if $\nx_- =\rm$,\quad $\bv.\bn_- =\frac{\rm}{2}$
\item
if $\nx_+ =\rp$,\quad $\bv.\bn_+ =\frac{\rm}{2}$
\item
for $i$ s.t. $|x_i-x_{i+1}|=\nx_-$,\quad $v_i-v_{i+1}=x_i-x_{i+1}$ 
hence $\bv.\bn_{i-}(\bx)=\frac{\nx_- - \nv_-}{\sqrt{3}} \ge \frac{\rm - \frac{\rm}{2}}{\sqrt{3}} = \frac{\rm}{2\sqrt{3}}$
\item
for $i$ s.t. $|x_i-x_{i+1}|=\nx_+$,\quad $v_i-v_{i+1}=x_{i+1}-x_i$ 
hence $\bv.\bn_{i+}(\bx)=\frac{\nx_+ + \nv_+}{\sqrt{3}} \ge \frac{\rm - \frac{\rm}{2}}{\sqrt{3}} = \frac{\rm}{2\sqrt{3}}$
\end{itemize}
We now proceed with the case  $\nx_- = \nx_+ < \rp$. Since $|x_i-x_{i+1}|=\nx_-=\nx_+$ and $v_i-v_{i+1}=x_i-x_{i+1}$ for each $i$, we have~:
\begin{itemize}
\item
if $\nx_- =\rm$,\quad $\bv.\bn_- =\frac{\rm}{2}$
\item
$\disp \bv.\bn_= = \frac{\frac{3}{2}\rp-\frac{\rm}{2}}{\sqrt{2}}  \ge \frac{\rp}{\sqrt{2}}$
\item
$\disp \bv.\bn_{i-}(\bx) =\frac{\nx_- - \nv_-}{\sqrt{3}} \ge \frac{\rm - \frac{\rm}{2}}{\sqrt{3}} = \frac{\rm}{2\sqrt{3}}$ and 
$\disp \bv.\bn_{i+}(\bx)=\frac{-\nx_- + \nv_+}{\sqrt{3}} \ge \frac{-\rp + \frac{3}{2}\rp}{\sqrt{3}} \ge \frac{\rp}{2\sqrt{3}}$
\end{itemize}
In the last case $\nx_- = \nx_+ = \rp$, one has $|x_i-x_{i+1}|=\rp$ and $v_i-v_{i+1}=x_{i+1}-x_i$ for each $i$~:
\begin{itemize}
\item
$\bv.\bn_+ =\frac{\rp}{2}$
\item
$\bv.\bn_= = \frac{\frac{-\rp}{2} + \frac{3}{2}\rp}{\sqrt{2}} \ge \frac{\rp}{\sqrt{2}}$
\item 
$\disp \bv.\bn_{i-}(\bx)=\frac{-\rp + \frac{3}{2}\rp}{\sqrt{3}} \ge \frac{\rp}{2\sqrt{3}}$ and 
$\disp \bv.\bn_{i+}(\bx)=\frac{\rp - \frac{\rp}{2}}{\sqrt{3}}=\frac{\rp}{2\sqrt{3}}$
\end{itemize}

So all these scalar products are larger than $\frac{\rm}{2\sqrt{3}}$.
As a consequence, assumption \emph{(iv)} in proposition \ref{propcritere} is satisfied with 
$\disp \beta_0 = \frac{\rm}{\rp n \sqrt{6n}}$.
Choosing $\beta_{i-}=\beta_{i+}=\beta_0^2/4$ hence $\delta_{i-}=\delta_{i+}=\frac{\rm\sqrt{3}}{64}\beta_0^4$
we obtain, thanks to proposition \ref{propcritere}, that $\A_c$ has the Uniform Exterior Sphere property with constant 
$\alpha_{\A_c}=\frac{\rm^2}{2\rp n \sqrt{2n}}$ and the Uniform Normal Cone property with constants 
$\delta_{\A_c}=\frac{\rm^5}{2^9 3 \sqrt{3}\rp^4 n^6}$ and $\beta_{\A_c}=\sqrt{1-\frac{\rm^2}{24\rp^2 n^3}}$
\end{proof2}

{\em Acknowledgments~:} The author thanks Sylvie Roelly for interesting discussions and helpful comments.


\begin{thebibliography}{99}

\bibitem{Fenchel}
W. Fenchel
Convex cones, sets, and functions.
Lecture notes, Princeton Univ. (1953).

\bibitem{FR2}
M. Fradon and S. Roelly, 
{\em Infinite system of Brownian balls with interaction: the non-reversible case},
ESAIM: P\&S 11 (2007) 55-79 

\bibitem{FR3}
M. Fradon and S. R\oe lly,
{\em Infinite system of Brownian balls : Equilibrium measures are canonical Gibbs}, 
Stochastics and Dynamics, Vol. 6 No. 1 (2006) 97-122.

\bibitem{IkedaWatanabe}
N. Ikeda and S. Watanabe,
Stochastic Differential Equations and Diffusion Processes.
North Holland - Kodansha (1981)

\bibitem{LionsSznitman}
P. L. Lions and A. S. Sznitman,
{\em Stochastic Differential Equations with Reflecting Boundary Conditions}, 
Com. Pure and Applied Mathematics {\bf  37} (1984) 511-537

\bibitem{SaishoTanakaBrownianBalls}      
Y. Saisho and H. Tanaka,
{\em Stochastic Differential Equations for Mutually Reflecting Brownian Balls},
Osaka J. Math. {\bf 23} (1986) 725-740

\bibitem{SaishoTanakaSymmetry}      
Y. Saisho and H. Tanaka,
{\em On the symmetry of a Reflecting Brownian Motion Defined by Skorohod's Equation for a Multi-Dimensional domain},
Tokyo J. Math. {\bf 10} (1987) 419-435

\bibitem{SaishoSolEDS}
Y. Saisho,
{\em Stochastic Differential Equations for multi-dimensional domain with reflecting boundary}, 
Probability Theory and Related Fields {\bf 74} (1987) 455-477.

\bibitem{SaishoMolecules2}
Y. Saisho,
{\em A model of the random motion of mutually reflecting molecules in $\R^d$}, 
Kumamoto Journal of Mathematics {\bf 7} (1994) 95-123.

\bibitem{Skorokhod1et2}
A. V. Skorokhod
{\em Stochastic equations for diffusion processes in a bounded region 1, 2}, 
Theor. Veroyatnost. i Primenen.  {\bf 6} (1961), 264-274;  {\bf 7} (1962) 3-23.

\bibitem{Tanaka}
H. Tanaka,
{\em Stochastic Differential Equations with Reflecting Boundary Conditions in Convex Regions}
Hiroshima Math. J., {\bf 9}, (1979), 163-177.


\end{thebibliography}
\end{document}